\documentclass{amsart}
\usepackage{amsmath, amssymb}
\usepackage{url, tikz-cd}
\usepackage[mathscr]{euscript}

\newtheorem*{theorem*}{Theorem}
  
\newtheorem{theorem}{Theorem}[section]
\newtheorem{lemma}[theorem]{Lemma}

\newtheorem*{corollary*}{Corollary}
\newtheorem{fact}[theorem]{Fact}

\newtheorem{claim}[theorem]{Claim}

\theoremstyle{definition}

\newtheorem{remark}[theorem]{Remark}
\newtheorem*{remark*}{Remark}
\newtheorem{definition}[theorem]{Definition}
\newtheorem*{definition*}{Definition}

\def \U {\mathcal U}
\def \C {\mathcal C}
\def \PP {\mathbb P}
\def \AA {\mathbb A}

\def\PGL{\operatorname{PGL}}
\def\DCF{\operatorname{DCF}}

\def\ccm{\operatorname{CCM}}

\def\dcl{\operatorname{dcl}}
\def\acl{\operatorname{acl}}

\def\tp{\operatorname{tp}}

\def\aut{\operatorname{Aut}}

\def\C{\mathcal C}

\def\nmdeg{\operatorname{nmdeg}}

\def\Ind#1#2{#1\setbox0=\hbox{$#1x$}\kern\wd0\hbox to 0pt{\hss$#1\mid$\hss}
\lower.9\ht0\hbox to 0pt{\hss$#1\smile$\hss}\kern\wd0}
\def\ind{\mathop{\mathpalette\Ind{}}}
\def\Notind#1#2{#1\setbox0=\hbox{$#1x$}\kern\wd0\hbox to 0pt{\mathchardef
\nn=12854\hss$#1\nn$\kern1.4\wd0\hss}\hbox to
0pt{\hss$#1\mid$\hss}\lower.9\ht0 \hbox to
0pt{\hss$#1\smile$\hss}\kern\wd0}
\def\nind{\mathop{\mathpalette\Notind{}}}

\title{The degree of nonminimality is at most~$2$}

\author{James Freitag}
\address{James Freitag\\
University of Illinois Chicago\\ 
Department of Mathematics, Statistics,
and Computer Science\\ 
851 S. Morgan Street\\
Chicago, IL, 60607-7045\\
USA}
\email{jfreitag@uic.edu}

\author{R\'emi Jaoui}
\address{R\'emi Jaoui\\
Albert-Ludwigs Universität Freiburg\\
Abteilung für Mathematische Logik, Mathematisches Institut\\ Ernst-Zermelo-Straße 1, D-79104 Freiburg\\ Germany.}
\email{remi.jaoui@math.uni-freiburg.de}

\author{Rahim Moosa}
\address{Rahim Moosa\\
University of Waterloo\\
Department of Pure Mathematics\\
200 University Avenue West\\
Waterloo, Ontario \  N2L 3G1\\
Canada}
\email{rmoosa@uwaterloo.ca}


\thanks{The first author was partially supported by NSF grant DMS-1700095 and NSF CAREER award 1945251.
The second author was partially supported by the ANR-DFG program GeoMod (Project number 2100310201).
The third author was partially supported by an NSERC DG}

\date{\today}

\begin{document}

\begin{abstract}
It is shown that if $p\in S(A)$ is a complete type of Lascar rank at least~$2$, in the theory of differentially closed fields of characteristic zero, then there exists a pair of realisations $a_1, a_2$ such that $p$ has a nonalgebraic forking extension over $Aa_1a_2$.
Moreover, if $A$ is contained in the field of constants then~$p$ already has a nonalgebraic forking extension over $Aa_1$.
The results are also formulated in a more general setting.
\end{abstract}

\maketitle

\section{Introduction}

\noindent
In~\cite{nmdeg}, motivated by the search for general techniques that might aid in proving strong minimality for certain algebraic differential equations,  the first and third authors introduced {\em degree of nonminimality} as a measure of how many parameters are needed to witness that a type is {\em not} minimal.
Working  in a sufficiently saturated model of a stable theory eliminating imaginaries, here is a precise formulation:

\begin{definition}
Suppose $p\in S(A)$ is a stationary type with $U(p)>1$.
The {\em degree of nonminimality of $p$}, denoted by $\nmdeg(p)$, is the least positive integer $d$ such that there exist realisations $a_1,\dots, a_d$ of $p$ and a nonalgebraic forking extension of~$p$ over $Aa_1,\dots,a_d$.
If $U(p)\leq 1$ then we set $\nmdeg(p)=0$ by convention.
\end{definition}

Using an analysis of the multiple transitivity of binding group actions, it was shown in~\cite{nmdeg} that $\nmdeg(p)\leq U(p)+1$ in the theory of differentially closed fields of characteristic zero ($\DCF_0$).
Bounds on the degree of nonminimality have played a significant role in recent proofs of strong minimality; of the generic  differential equation in~\cite{freitag-devilbiss} and of the differential equations satisfied by the Schwarz triangle functions in~\cite{cdfn}.
Based on a maturing of the techniques used in~\cite{nmdeg}, and informed by the approach taken in~\cite{c3c2} to a related problem, we give in this note a short proof of a dramatic improvement to that bound:

\begin{theorem*}
Suppose $T=\DCF_0$ and $p$ is a complete stationary type of finite rank.
Then $\nmdeg(p)\leq 2$.
Moreover, if $p$ is over constant parameters then $\nmdeg(p)\leq 1$.
\end{theorem*}

The bound is sharp; see~\cite[Example~4.2]{nmdeg} for types of nonminimality degree~$2$.

The argument we give for the main clause, namely that $\nmdeg(p)\leq 2$, works equally well in $\DCF_{0,m}$, the theory of differentially closed fields in $m$ commuting derivations, and in $\ccm$, the theory of compact complex manifolds.
All one needs is that $T$ be totally transcendental, eliminate imaginaries, eliminate the ``there exists infinitely many" quantifier, and admit a $0$-definable pure algebraically closed field to which every non locally modular minimal type is nonorthogonal.
In $\DCF_{0,m}$ that pure algebraically closed field is the field of constants and in $\ccm$ it is the (interpretation in~$\U$ of the) complex field living on the projective line.

The ``moreover" clause of the theorem, however, does make use of the fact that, in $\DCF_0$, the binding group of a type over the constants and internal to the constants cannot be centerless.

The most general setting for the results is articulated, for the record, in Section~\ref{general}.

\begin{remark}
A corollary of our theorem is a significant improvement to the main result of~\cite{freitag-devilbiss}, where it was shown that generic algebraic differential equations of order $h\geq 2$ and degree at least $2(h + 2)$ are strongly minimal.
The proof in~\cite{freitag-devilbiss} used that $\nmdeg(p)\leq U(p)+1$.
The same proof, but using the improved bound of  $\nmdeg(p)\leq 2$ obtained here, allows one to replace $2(h + 2)$ by~$6$ in that result.
\end{remark}

\bigskip
\section{The proof}

\noindent
We work in a fixed sufficiently saturated model $\U$ of a complete totally transcendental theory $T$ eliminating imaginaries and the ``there exists infinitely many" quantifier, with $\C$ a $0$-definable pure algebraically closed field such that every non locally modular minimal type is nonorthogonal to $\C$.

Maybe the first thing to observe is that the degree of nonminimality is invariant under interalgebraicity.
Here we use the following, possibly nonstandard but unambigious, terminology:

\begin{definition}
Complete types $p,q\in S(A)$ are said to be {\em interalgebraic}
if for each (equivalently some) $a\models p$ there exists $b\models q$ such that $\acl(Aa)=\acl(Ab)$
.
\end{definition}

That $\nmdeg(p)=\nmdeg(q)$ when $p$ and $q$ are interalgebraic is more or less immediate from the definitions; see for example~\cite[Lemma~3.1(c)]{nmdeg}.

The following consequences of $\nmdeg>1$ were observed in~\cite{nmdeg}, but we include some details here for the sake of completeness:

\begin{fact}
\label{nmdeg>1}
Suppose $p\in S(A)$ is stationary of finite rank with $\nmdeg(p)>1$.
Then $p$ is interalgebraic with a stationary type $q\in S(A)$ such that $q$ is $\C$-internal and $q^{(2)}$ is weakly $\C$-orthogonal.
\end{fact}

\begin{proof}
Note, first of all, that 
\begin{itemize}
\item[($*$)]
if $a\models p$ and $b\in\acl(Aa)\setminus\acl(A)$ then $a\in\acl(Ab)$.
\end{itemize}
Indeed, if $a'$ realises the nonforking extension of $p$ to $Aab$ then $\tp(a'/Aa)$ is a forking extension of $p$.
Since $\nmdeg(p)>1$ we must have that $a'\in\acl(Aa)$, from which it follows that $a'\in\acl(Ab)$, and hence $a\in\acl(Ab)$.

In the finite rank setting, condition~($*$), which is a weak form of exchange, implies that either $p$ is interalgebraic with a  locally modular minimal type, or $p$ is almost internal to a non locally modular minimal type -- see~\cite[Proposition~2.3]{moosa2014some}.
The former is impossible as $U(p)>1$, and by assumption on $T$ the latter implies $p$ is almost $\C$-internal.
We thus find a stationary $\C$-internal $q\in S(A)$ that is interalgebraic with~$p$.
Note that $\nmdeg(q)>1$ as well.

Suppose that $q$ is not weakly $\C$-orthogonal.
Since the induced structure on $\C$, namely that of a pure algebraically closed field, eliminates imaginaries, this failure of weak $\C$-orthogonality will be witnessed by some $b\models q$ and $c\in \C$ such that $c\in\dcl(Ab)\setminus \acl(A)$.
By~($*$) applied to $q$ this would force $b\in\acl(Ac)$, contradicting $U(q)>1$.
So $q$ is weakly $\C$-orthogonal.
In particular, as it is $\C$-internal, $q$ is isolated.
We let $\Omega$ be the definable set $q(\U)$.

Now suppose that $q^{(2)}$ is not weakly $\C$-orthogonal.
Then there are independent $b_1,b_2$ realising $q$ and $c\in \C$ such that $c\in\dcl(Ab_1b_2)\setminus\acl(A)$.
Note that $b_2\notin\acl(Ab_1c)$ as $U(b_2/Ab_1)=U(q)>1$.
So there is a partial $Ab_1$-definable function $f:\Omega\to\C$ with infinite image and infinite generic fibre.
It follows, by elimination of the ``there exists infinitely many" quantifier, that all but finitely many of the fibres are infinite.
As $\C\cap\acl(A)$ is infinite (it is an algebraically closed subfield of $\C$), there exists $b\in \Omega\setminus\acl(Ab_1)$ such that $f(b)\in\acl(A)$.
If $b\ind_Ab_1$ then $\tp(b/Ab_1)=\tp(b_2/Ab_1)$ contradicting the fact that $f(b_2)=c\notin\acl(A)$.
So $b\nind_Ab_1$.
That is, $\tp(b/Ab_1)$ is a nonalgebraic forking extension of $q$.
But this contradicts $\nmdeg(q)>1$.
Hence $q^{(2)}$ is weakly $\C$-orthogonal.
\end{proof}

The following improvement to Fact~\ref{nmdeg>1} was {\em not} remarked in~\cite{nmdeg}.

\begin{lemma}
\label{lem:nmdeg>1}
Suppose $p\in S(A)$ is stationary of finite rank with $\nmdeg(p)>1$.
Then $p$ is interalgebraic with some stationary $q\in S(A)$ such that
\begin{itemize}
\item[(a)]
$q$ is $\C$-internal,
\item[(b)]
$q^{(2)}$ is weakly $\C$-orthogonal, and,
\item[(c)]
any two distinct realisations of $q$ are independent over $A$.
\end{itemize}
\end{lemma}

\begin{proof}
Suppose $a,b$ are realisations of $p$ such that $a\nind_Ab$.
If $a\notin\acl(Ab)$ then $\tp(a/Ab)$ is a nonalgebraic forking extension of $p$, contradicting $\nmdeg(p)>1$.
Similarly, we must have $b\in\acl(Aa)$.
In other words, $a\nind_Ab$ if and only if $\acl(Aa)=\acl(Ab)$.
In particular, $a\nind_Ab$ is an equivalence relation on $p(\U)$, which we now denote by $E$.

Applying Fact~\ref{nmdeg>1}, we may assume that $p$ is $\C$-internal and $p^{(2)}$ is weakly $\C$-orthogonal.
In particular, both $p$ and $p^{(2)}$ are isolated, say by the $L_A$-formulae $\phi(x)$ and $\psi(x,y)$, respectively.
Note then, that $\phi(x)\wedge\phi(y)\wedge \neg\psi(x,y)$ defines the forking relation $E$.
So $E$ is an $A$-definable equivalence relation.

Each class of $E$ is finite.
Indeed, if $a\models p$ has an infinite $E$-class then there is $b\in p(\U)\setminus\acl(Aa)$ with $aEb$.
But that means that $\tp(b/Aa)$ is a nonalgebraic forking extension of $p$, contradicting $\nmdeg(p)>1$.

Fixing $a\models p$, let $e:=a/E$ and $q:=\tp(e/A)$.
Note that $e\in\dcl(Aa)$, and so we still have that $q$ is $\C$-internal and $q^{(2)}$ is weakly $\C$-orthogonal.
Also, as the $E$-classes are finite, $p$ and $q$ are interalgebraic.
So it remains to show that any two distinct realisations of $q$ are independent.
Suppose $e'\models q$ with $e'\neq e$.
Then $e'=a'/E$ for some $a'\models p$ such that $\neg(aE a')$.
That is $a\ind_Aa'$.
As $\acl(Aa)=\acl(Ae)$ and $\acl(Aa')=\acl(Ae')$, we have that $e\ind_Ae'$, as desired.
\end{proof}

We now work toward a proof of the main clause of the Theorem.
That is, fixing a finite rank stationary type $p\in S(A)$, we wish to show that $\nmdeg(p)\leq 2$.
Let $\overline p$ denote the unique extension of $p$ to $\acl(A)$.
It is immediate from the definition that $\nmdeg(\overline p)=\nmdeg(p)$.
We may therefore assume that $A=\acl(A)$.
Let $k:=A\cap\C$, it is an algebraically closed subfield of $\C$.

In order to prove that $\nmdeg(p)\leq 2$ we may of course assume that $\nmdeg(p)>1$.
Hence, by Lemma~\ref{lem:nmdeg>1}, we can further reduce to the case that $p$ is $\C$-internal, $p^{(2)}$ is weakly $\C$-orthogonal, and any two distinct realisations of $p$ are independent over~$A$.

Let $\Omega:=p(\U)$ and let $G:=\aut(p/\C)$ be the binding group of $p$ relative to~$\C$.
So $(G,\Omega)$ is an $A$-definable faithful group action.
The action is transitive because $p$ is weakly $\C$-orthogonal.
Weak $\C$-orthogonality of $p$ also implies, along with $A=\acl(A)$, that $G$ is connected.
The fact that $p^{(2)}$ is weakly $\C$-orthogonal implies that $G$ acts transitively on $p^{(2)}(\U)$.
But $p^{(2)}(\U)=\Omega^2\setminus\Delta$ where $\Delta$ is the diagonal, because any two distinct realisations of $p$ are independent over $A$.
So $(G,\Omega)$ is a $2$-transitive connected $A$-definable homogeneous space.

Now, the binding group action of any $\C$-internal type is isomorphic to the $\C$-points of an algebraic group action, though possibly over additional parameters.
More precisely, let $M\preceq\U$ be a prime model over $A$.
Note that $M\cap\C=k$.
There exists an algebraic homogeneous space $(\overline G,\overline \Omega)$ defined over $k$, and an $M$-definable isomorphism $\alpha:(G,\Omega)\to(\overline G(\C),\overline \Omega(\C))$.

In particular, $(\overline G,\overline \Omega)$ is a $2$-transitive connected algebraic homogeneous space.
This is a very restrictive condition; a theorem of Knop~\cite{knop1983mehrfach} tells us that $(\overline G,\overline \Omega)$ is either isomorphic to the action of $\PGL_{n+1}$ on $\PP^n$, or is isomorphic to the action of an algebraic subgroup of the group of affine transformations on $\AA^n$, for some $n>1$.
In either case we have a notion of {\em collinearity} which is preserved by the group action.
That is, given distinct $u,v\in\overline \Omega(\C)$ we can talk about the line $\ell_{u,v}\subseteq\overline \Omega(\C)$ connecting $u$ and $v$, and for all $g\in \overline G(\C)$ we have that $g\ell_{u,v}=\ell_{gu,gv}$.

Fix distinct $a,b\in \Omega$, and consider the set $X:=\alpha^{-1}(\ell_{\alpha(a),\alpha(b)})$.
Then $X$ is a rank~$1$ $Mab$-definable subset of $\Omega$.

\begin{claim}
There is a finite tuple $c$ from $\C$ such that $X$ is $Aabc$-definable.
\end{claim}

\begin{proof}
It suffices to show that if $\sigma\in\aut_{Aab}(\U/\C)$, that is, if $\sigma$ is an automorphism of $\U$ that fixes $A\cup\{a,b\}\cup\C$ point-wise, then $\sigma(X)=X$.
Now, the restriction of $\sigma$ to $\Omega$ is an element of the binding group, say $g_{\sigma}\in G$, which fixes $a$ and $b$.
Hence $\alpha(g_{\sigma})\in \overline G(\C)$ fixes $\alpha(a)$ and $\alpha(b)$, and hence preserves the line $\ell_{\alpha(a),\alpha(b)}$.
It follows that
\begin{eqnarray*}
\alpha(\sigma(X))
&=&
\alpha\big(g_\sigma(\alpha^{-1}(\ell_{\alpha(a),\alpha(b)}))\big)\\
&=&
\alpha(g_\sigma)(\ell_{\alpha(a),\alpha(b)})\\
&=&
\ell_{\alpha(a),\alpha(b)}.
\end{eqnarray*}
Applying $\alpha^{-1}$ to both sides we obtain that $\sigma(X)=X$, as desired.
\end{proof}

Let $\theta(x,y)$ be an $L_{Aab}$-formula such that $X=\theta(\U,c)$.
If, in addition, we chose $a,b\in \Omega(M)$, then $X$ and $\theta(x,y)$ are over $M$, and it follows that there is $c'\in M\cap\C$ such that $X=\theta(\U,c')$.
But $M\cap\C=k\subseteq A$, so  that this witnesses the definability of $X$ over $Aab$.

We have thus found $a,b\in \Omega$ and an $Aab$-definable subset $X\subseteq \Omega$ of rank~$1$.
Since $U(p)>1$, the generic type of $X$ over $Aab$ is a nonalgebraic forking extension of $p$.
Since $a$ and $b$ realise $p$, this witnesses that $\nmdeg(p)=2$.

This completes the proof of the main clause of the Theorem.

For the ``moreover" clause, we return to the particular setting of $T=\DCF_0$ and~$\C$ the field of constants.
We make the additional assumption that $A\subseteq\C$ and show that $\nmdeg(p)>1$ leads to a contradiction.
Indeed, that $(G,\Omega)$ is $2$-transitive forces $G$ to be centerless; see for example the elementary argument at the beginning of the proof of Satz~2 in~\cite{knop1983mehrfach}.
But, in $\DCF_0$, the binding group of a type that is $\C$-internal and over constant parameters cannot be centerless; see for example the proof of Theorem~3.9 in~\cite{c3c2}.
This contradiction proves that $\nmdeg(p)\leq 1$.
\qed

\bigskip
\section{Some remarks on the assumptions}
\label{general}

\noindent
We carried out the above proof under assumptions on $T$ that were suitable for generalisation to both $\DCF_{0,m}$ and $\ccm$.
But it may be worth recording the minimal hypotheses on $T$ required for the proofs to go through.
We leave it to the reader to inspect those proofs and verify that what is actually proved are the following two statements:

\begin{theorem}
Suppose $T$ is a complete totally transcendental theory eliminating imaginaries and the ``there exists infinitely many" quantifier.
Let $\U\models T$ be a sufficiently saturated model and $A\subseteq \U$ a parameter set.
\begin{itemize}
\item[(a)]
Suppose each non locally modular minimal type in $T$ is nonorthogonal to some $A$-definable pure algebraically closed field.
Then $\nmdeg(p)\leq 2$ for all stationary $p\in S(A)$ of finite rank.
\item[(b)]
Suppose there exists a collection $\{\C_i:i\in I\}$ of $A$-definable non locally modular strongly minimal sets such that  each non locally modular minimal type in $T$ is nonorthogonal to $\C_i$ for some $i\in I$, and such that for all $i\in I$,
\begin{itemize}
\item[(i)]
 $\C_i\cap\acl(A)$ is infinite, and,
 \item[(ii)]
 for all  $\C_i$-internal $q\in S(A)$, the binding group $\aut(q/\C_i)$ has a nontrivial center.
 \end{itemize}
 Then $\nmdeg(p)\leq 1$ for all stationary $p\in S(A)$ of finite rank.
 \end{itemize}
\end{theorem}

\end{document}